\documentclass[12pt,a4paper,fleqn]{article}
\usepackage{a4wide,amsfonts,amsmath,latexsym,amssymb,euscript,graphicx,units,mathrsfs}

\usepackage{graphicx}
\usepackage{color}
\usepackage{amssymb}
\usepackage[T1]{fontenc}
\usepackage{latexsym}
\usepackage{xypic}
\usepackage{eufrak}
\usepackage{euscript}
\usepackage{amsfonts,amsmath}
\usepackage{verbatim}
\usepackage{fancyhdr}
\usepackage[english]{babel}
\usepackage{mathrsfs}
\usepackage{units}

\newtheorem{prop}{Proposition}[section]
\newtheorem{cor}[prop]{Corollary}

\newtheorem{thm}[prop]{Theorem}

\renewcommand{\geq}{\geqslant}
\def\leq{\leqslant}

\newcommand{\N}{\mathbb{N}}

\newcommand{\R}{\mathbb{R}}

\def\HH{\EuFrak H}

\def\e{\varepsilon}

\def\1{{\mathbf{1}}}

\def\1{{\mathbf{1}}}
\def\0.5{{\frac{1}{2}}}

\newcommand{\qed}{\nopagebreak\hspace*{\fill}
{\vrule width6pt height6ptdepth0pt}\par}

\newcounter{rea}
\begin{document}

\title{Absolute continuity and convergence of  densities for random vectors on Wiener chaos}
\author{Ivan Nourdin\thanks{
Universit\'e de Lorraine, Institut \'Elie Cartan de Lorraine, UMR 7502, Vandoeuvre-l\`es-Nancy, F-54506, France and CNRS, Institut \'Elie Cartan de Lorraine, UMR 7502, Vandoeuvre-l\`es-Nancy, F-54506, France, {\tt inourdin@gmail.com}; IN is partially supported by the
ANR grants ANR-09-BLAN-0114 and ANR-10-BLAN-0121.}, \, David Nualart\thanks{
Department of Mathematics,
University of Kansas,
Lawrence, Kansas, 66045 USA, {\tt nualart@math.ku.edu}; DN is supported by the NSF grant DMS-1208625.}  \, and \, Guillaume Poly\thanks{Laboratoire d'Analyse et de Math\'ematiques Appliqu\'ees, UMR 8050, Universit\'e Paris-Est Marne-
la-Vall\'ee, 5 Bld Descartes, Champs-sur-Marne, 77454 Marne-la-Vall\'ee Cedex 2, France, {\tt guillaume.poly@crans.org}.}}
\maketitle

\begin{abstract}
The aim of this paper is to establish some new results on the absolute continuity and the convergence in total variation for a sequence of $d$-dimensional vectors whose components belong to a finite sum of Wiener chaoses. First we   show that the probability that the determinant of the Malliavin matrix of such vectors  vanishes is zero or one, and this probability equals to one is equivalent to say that the vector takes values in the set of zeros of a polynomial. We provide a bound for the degree of this annihilating polynomial  improving a result by Kusuoka  \cite{kusuoka}. On the other hand, we show that the convergence in law implies the convergence in total variation, extending to the multivariate case a recent result by Nourdin and Poly \cite{NouPol}. This follows from an inequality relating the total variation distance with the Fortet-Mourier distance. Finally, applications to some particular cases are discussed.

\end{abstract}

\section{Introduction}
The purpose of this paper is to establish some new results on the absolute continuity and the convergence of the densities in some  $L^p(\mathbb{R}^d)$  for a sequence of $d$-dimensional random vectors whose components belong to finite sum of Wiener chaos.
These result generalize previous works by Kusuoka \cite{kusuoka} and by Nourdin and Poly
\cite{NouPol}, and are based on a combination of the techniques of Malliavin calculus, the Carbery-Wright inequality and some recent work on algebraic dependence for a family of polynomials.

Let us describe our main results. Given two $d$-dimensional random vectors $F$ and $G$, we denote by $d_{TV}(F,G)$  the total variation distance between the laws of $F$ and $G$, defined by
\[
d_{TV}(F,G) =\sup_{A \in \mathcal{B}(\mathbb{R}^d)} | P(F\in A)-P(G\in A)|,
\]
where the supremum is taken over all Borel sets $A$  of $\mathbb{R}^d$.
There is an equivalent formulation for $d_{TV}$, which is often useful:
\[
d_{TV}(F,G) =\frac12 \,\sup_\phi |E[\phi(F)] - E[\phi(G)] |,
\]
where the supremum is taken over all measurable functions $\phi:\mathbb{R}^d\rightarrow \mathbb{R}$ which are bounded by $1$.
It is also well-known (Scheff\'e's Theorem) that, when $F$ and $G$ both have a law which is absolutely continuous with respect to the Lebesgue measure on $\R^d$, then
\[
d_{TV}(F,G) =\frac12 \int_{\R^d} |f(x)-g(x)|dx,
\]
with $f$ and $g$ the densities of $F$ and $G$ respectively.
On the other hand,
we denote by $d_{FM}(F,G)$ the Fortet-Mourier distance, given by
\[
d_{FM}(F,G)= \sup_\phi |E[\phi(F)] - E[\phi(G)] |,
\]
where the supremum is taken over all $1$-Lipschitz functions $\phi:\mathbb{R}^d\rightarrow \mathbb{R}$ which are bounded by $1$. It is well-known that $d_{FM}$ metrizes the convergence in distribution.

Consider a sequence of random vectors $F_n =(F_{1,n}, \dots, F_{d,n})$ whose components belong to $\oplus_{k=0}^q \mathcal{H}_k$, where $\mathcal{H}_k$ stands for the $k$th Wiener chaos, and assume that $F_n$ converges in distribution towards a random variable $F_{\infty}$. Denote by $\Gamma(F_n)$ the Malliavin matrix of $F_n$, and assume that $E[\det\Gamma(F_n)]$ is bounded away from zero.
 Then we prove that there exist  constants $c,\gamma>0$ (depending on $d$ and $q$) such that, for any $n\geq 1$,
 \begin{equation} \label{hh1}
 d_{TV}(F_n ,F_\infty) \leq c d_{FM} (F_n, F_\infty) ^\gamma.
 \end{equation}
 So, our result implies that the sequence $F_n$ converges not only in law but also in total variation.
 In \cite{NouPol} this result has been proved for $d=1$. In this case $\gamma=\frac 1{2q+1}$, and one only needs that $F_\infty$ is not identically zero, which turns out to be equivalent to the fact that the law of $F_\infty$ is absolutely continuous. This equivalence is not true for $d\geq 2$. The proof of this result is based on the Carbery-Wright inequality for the law of a polynomial on Gaussian random variables and also on the integration-by-parts formula of Malliavin calculus. In the multidimensional case we make use of the integration-by-parts formula based on the Poisson kernel developed by Bally and Caramelino in  \cite{BaCa1}.

The convergence in total variation is very strong, and should
not be expected from the mere convergence in law without some additional structure. For instance, there is a celebrated theorem of
Ibragimov (see, e.g., Reiss \cite{reiss}) according to which, if $F_n,F_\infty$ are continuous random variables with densities $f_n,f_\infty$ that are {\it unimodal}, then $F_n\to F_\infty$ in law if and only if $ d_{TV}(F_n ,F_\infty)\to 0$.
Somehow, our inequality (\ref{hh1}) thus appears as unexpected.
Several consequences are detailed in Section 5.
Furthermore, bearing in mind that the convergence in total variation is equivalent to the convergence of the densities in $L^1(\mathbb{R}^d)$, we improve this results by proving that under the above assumptions on the sequence $F_n$, the densities converge in $L^p(\mathbb{R}^d)$ for some explicit $p>1$ depending solely on $d$ and $q$.

 Motivated by the above inequality (\ref{hh1}), in the first part of the paper we discuss the absolute continuity of the law of a $d$-dimensional random vector $F=(F_1, \dots, F_d)$ whose components belong to a  finite sum of Wiener chaoses $\oplus_{k=1}^q \mathcal{H}_k$.
 Our main result says that the three following conditions are equivalent:
 \begin{enumerate}
 \item The law of $F$ is not absolutely continuous with respect to the Lebesgue measure on $\mathbb{R}^d$.
 \item  There exists a nonzero polynomial $H$  in $d$ variables of degree at most $dq^{d-1}$   such that $H(F)=0$.
 \item $E[\det \Gamma(F)]=0$.
\end{enumerate}
Notice that the criterion  of the Malliavin calculus for the absolute continuity  of the law of a random vector $F$  says that $\det \Gamma(F)>0$ almost surely implies the absolute continuity of the law of $F$.  We prove the stronger result that $P(\det \Gamma(F) =0)$ is zero or one; as a consequence, $P(\det \Gamma(F) >0)=1$ turns out to be equivalent to the absolute continuity. The equivalence with condition 2 improves a classical result by Kusuoka (\cite{kusuoka}), in the sense that we provide a simple proof of the existence of the annihilating polynomial based on a recent result by Kayal \cite{kayal} and  we give an upper bound for the degree of this polynomial. Also, it is worthwhile noting that, compared to condition 2, condition 3 is often easier to check in practical situations, see also the end of Section \ref{sec-ac}.

The paper is organized as follows. Section 2 contains some preliminary material on Malliavin calculus, the Carbery-Wright inequality and the  results on algebraic dependence that will be used in the paper. In Section 3 we provide equivalent conditions for absolute continuity in the case of a random vector in a sum of Wiener chaoses. Section 4 is devoted to establish the inequality  (\ref{hh1}), and also the convergence in $L^p(\mathbb{R}^d)$ for some $p$. Section 5 contains   applications of these results in some particular cases. Finally, we list two open questions in Section 6.

\section{Preliminaries}

This section contains some basic elements on Gaussian analysis that will be used throughout this paper. We refer the reader to the books  \cite{NouPec, nualartbook} for further details.

\subsection{Multiple stochastic integrals}

Let $\mathfrak{H}$  be a real separable Hilbert space. We denote by $X=\{X(h), h\in \mathfrak{H}\}$ an \textit{isonormal Gaussian process} over $\mathfrak{H}$. That means, $X$ is a centered Gaussian family of random variables defined in some probability space $(\Omega, \mathcal{F},P)$, with covariance given by
\[
E[X(h)X(g)]= \langle h,g \rangle_{\mathfrak{H}},
\]
for any $h,g \in \mathfrak{H}$. We also assume that $\mathcal{F}$ is generated by $X$.

For every $k\ge1$, we denote by $\mathcal{H}_k$ the $k$th \textit{Wiener chaos} of $X$ defined as the closed linear subspace of $L^2(\Omega)$ generated by the family of random variables $\{H_k(X(h)), h\in \mathfrak{H}, \|h\|_{\mathfrak{H}}=1\}$, where $H_k$ is the $k$th Hermite polynomial given by
\[
H_k(x)= (-1)^k e^{\frac {x^2}2} \frac {d^k} {dx^k} \left( e^{-\frac {x^2}2}\right).
\]
We write by convention $\mathcal{H}_0=\mathbb{R}$.
For any $k\geq 1$, we denote by $\mathfrak{H}^{\otimes k}$ the $k$th tensor product of $\mathfrak{H}$. Then, the mapping $I_k(h^{\otimes k})= H_k(X(h))$ can be extended to a linear isometry between the symmetric tensor product $\mathfrak{H}^{\odot k}$ (equipped with the modified norm $\sqrt{k!} \| \cdot \|_{\mathfrak{H}^{\otimes k}}$) and the $k$th Wiener chaos $\mathcal{H}_k$. For $k=0$ we write $I_0(x)=c$, $c\in \mathbb{R}$.  In the particular case $\mathfrak{H}=L^2(A,\mathcal{A}, \mu)$, where $\mu$ is a $\sigma$-finite measure without atoms, then $ \mathfrak{H}^{\odot k}$ coincides with the space $ L^2_s(\mu^k)$  of symmetric functions which are square integrable with respect to the product measure $\mu^k$, and for any $f\in \mathfrak{H}^{\odot k}$ the random variable $I_k(f)$ is the multiple stochastic integral of $f$ with respect to the centered Gaussian measure generated by $X$.

Any random variable $F\in L^2(\Omega)$ admits an orthogonal decomposition of the form
$F= \sum_{k=0}^\infty I_k(f_k)$, where $f_0=E[F]$, and the kernels $f_k\in  \mathfrak{H}^{\odot k}$ are uniquely determined by $F$.

Let $\{e_i, i\geq 1\}$ be a complete orthonormal system in $\mathfrak{H}$. Given $f\in \mathfrak{H}^{\odot k}$ and $g\in \mathfrak{H}^{\odot j}$, for every $r=0,\dots, k\wedge j$, the \textit{contraction} of $f$ and $g$ of order $r$ is the element of $\mathfrak{H}^{\otimes(k+j-2r)}$ defined by
\[
f\otimes_r g= \sum_{i_1, \dots, i_r=1}^\infty \langle f, e_{i_1} \otimes \cdots \otimes e_{i_r}
\rangle_{\mathfrak{H}^{\otimes r}} \otimes  \langle g, e_{i_1} \otimes \cdots \otimes e_{i_r}
\rangle_{\mathfrak{H}^{\otimes r}}.
\]
The contraction $f\otimes_r g$ is not necessarily symmetric, and we denote by $f \widetilde{\otimes}_rg$ its symmetrization.

\subsection{Malliavin calculus}
Let $\mathcal{S}$ be the set of all cylindrical random variables of the form
\[
F=g(X(h_1), \dots, X(h_n)),
\]
where $n\geq 1$, $h_i \in \mathfrak{H}$, and $g$ is infinitely differentiable such that all its partial derivatives have polynomial growth. The Malliavin derivative of $F$ is the element of $L^2(\Omega;\mathfrak{H})$ defined by
\[
DF= \sum_{i=1}^n \frac {\partial g}{\partial x_i}(X(h_1), \dots, X(h_n)) h_i.
\]
By iteration, for every $m\geq 2$, we define the $m$th derivative $D^mF$ which is an element of $L^2(\Omega; \mathfrak{H}^{\odot m})$. For $m\geq 1$ and $p\geq 1$, $\mathbb{D}^{m,p}$ denote the closure of $\mathcal{S}$ with respect to the norm $\| \cdot \|_{m,p}$ defined by
\[
\|F\|^p_{m,p} = E[|F|^p] + \sum_{j=1}^m E\left[ \|D^jF\|^p_{\mathfrak{H}^{\otimes j}}\right].
\]
We also set $\mathbb{D}^\infty = \cap_{m\geq 1} \cap_{p\geq 1} \mathbb{D}^{m,p}$.

As a consequence of the hypercontractivity property of the Ornstein-Uhlenbeck semigroup,
all the $\| \cdot \|_{m,p}$-norms are equivalent in a finite Wiener chaos. This is a basic result that will be used along the paper.

We denote by $\delta$ the adjoint of the operator $D$, also called the \textit{divergence operator}. An element $u\in L^2(\Omega;\mathfrak{H})$ belongs to the domain of $\delta$, denoted $\mathrm{Dom} \delta$, if $|E\langle DF,u \rangle_{\mathfrak{H}} | \leq c_u \|F\|_{L^2(\Omega)}$ for any $F\in \mathbb{D}^{1,2}$,  where $c_u$ is a constant depending only on $u$. Then, the random variable $\delta(u)$ is defined by the duality relationship
\begin{equation}\label{dual}
E[F\delta(u)]= E\langle DF,u \rangle_{\mathfrak{H}}.
\end{equation}
Given a random vector $F=(F_1, \dots, F_d)$ such that $F_i \in \mathbb{D}^{1,2}$, we denote  $\Gamma(F)$ the \textit{Malliavin matrix} of $F$, which is a random nonnegative definite matrix defined by
\[
\Gamma_{i,j}(F) = \langle DF_i, DF_j \rangle_{\mathfrak{H}}.
\]
If $F_{i,n}\in\mathbb{D}^{1,p}$ for some $p>1$ and any $i=1,\ldots,d$, and if $\det \Gamma (F)>0$ almost surely, then the law of $F$ is absolutely continuous with respect to the Lebesgue measure on $\mathbb{R}^d$ (see, for instance, \cite[Theorem 2.1.2]{nualartbook}). This is our basic criterion for absolute continuity in this paper.

\subsection{Carbery-Wright inequality}

Along the paper we will make use of the following inequality due to Carbery and Wright \cite{CW}: there is a universal constant $c>0$ such that, for any polynomial $Q:\mathbb{R}^n \rightarrow \mathbb{R}$ of degree at most $d$ and  any $\alpha>0$ we have
\begin{equation}  \label{ff1}
E[Q(X_1, \dots, X_n)^2] ^{\frac 1{2d}} P(|Q(X_1, \dots, X_n)| \leq \alpha) \leq cd\alpha^{\frac 1d},
\end{equation}
where $X_1, \dots, X_n$ are independent random variables with law $N(0,1)$.

\subsection{Algebraic dependence}

Let $\mathbb{F}$ be a field and $\mathbf{f}=(f_1, \dots, f_k) \in \mathbb{F}[x_1, \dots, x_n]$ be a set of $k$ polynomials of degree at most  $d$ in $n$ variables in the field $\mathbb{F}$. These polynomials are said to be \textit{algebraically dependent} if there exists a nonzero $k$-variate polynomial
$A(t_1, \dots, t_k) \in \mathbb{F}[t_1, \dots, t_k]$ such that $A(f_1, \dots, f_k)=0$. The polynomial $A$ is then called an $(f_1, \dots, f_k)$-\textit{annihilating polynomial}.

Denote by
\[
J_{\mathbf{f}} = \left( \frac {\partial f_i }{\partial x_j} \right) _{1\leq i \leq k, 1\leq j \leq n}
\]
the Jacobian matrix of the set of polynomials in $\mathbf{f}$. A classical result (see, e.g., Ehrenborg and Rota \cite{ER} for a proof) says that  $f_1, \dots, f_k$ are  algebraically independent if and only if the Jacobian matrix $J_{\mathbf{f}} $ has rank $k$.

Suppose that the polynomials $\mathbf{f}=(f_1, \dots, f_k) $ are algebraically dependent. Then the set of  $\mathbf{f}$-annihilating polynomials forms an ideal in the polynomial ring $\mathbb{F}[t_1, \dots, t_k]$.
In a recent work Kayal (see \cite{kayal}) has established some properties of this ideal. In particular (see \cite{kayal},  Lemma 7) he has proved that  if no proper subset of   $\mathbf{f}$ is algebraically  dependent, then the ideal of $\mathbf{f}$-annihilating polynomials  is generated by a single irreducible polynomial.  On the other hand (see  \cite{kayal}, Theorem 11) the degree of this generator is at most $k q^{k-1}$.

\section{Absolute continuity of the law of a system of multiple stochastic integrals}\label{sec-ac}

The purpose of this section is to extend a result by Kusuoka \cite{kusuoka} on the  characterization of the absolute continuity of a vector whose components are finite sums of multiple stochastic integrals, using techniques of Malliavin calculus. In what follows, the notation $\R[X_1,\ldots,X_d]$ stands for the set of $d$-variate polynomials over $\R$.

\begin{thm}\label{ac-vectormultiple}
Fix $q,d\geq 1$, and let $F=(F_1,\ldots,F_d)$  be a random vector such that $F_i\in\bigoplus_{k=1}^q \mathcal{H}_k$ for any $i=1,\ldots,d$. Let $\Gamma:=\Gamma(F)$ be the Malliavin matrix of $F$.
Then the following assertions are equivalent:
\begin{itemize}
\item[(a)] The law of $F$ is not absolutely continuous with respect to the Lebesgue measure on $\mathbb{R}^d$.
\item[(b)] There exists $H\in\R[X_1,\dots,X_d]\setminus\{0\}$ of degree at most $D=dq^{d-1}$ such that, almost surely,
\[
H(F_1,\dots,F_d)=0.
\]
\item[(c)] $E[\det\Gamma]=0$.
\end{itemize}
\end{thm}
{\it Proof of (a)$\Rightarrow$(c)}.
Let us prove $\neg \mbox(c) \Rightarrow \neg \mbox(a)$.
Set $N=2d(q-1)$
and let $\{e_k, k\geq 1\} $ be an orthonormal basis of $\mathfrak{ H}$. Since $\det \Gamma \in \bigoplus_{k=0}^N\mathcal{H}_k$,
there exists a sequence $\{Q_n, n\geq 1\}$ of real-valued polynomials of degree at most $N$ such that the random variables
$Q_n(I_1(e_1),\ldots,I_1(e_n))$ converge in $L^2(\Omega)$ and almost surely to $\det \Gamma$ as $n$ tends to infinity (see \cite[Theorem 3.1, first step of the proof]{NouPol} for an explicit construction).
Assume now that $E[\det \Gamma]>0$. Then for $n\geq n_0$,
$E[|Q_n(I_1(e_1),\ldots,I_1(e_n))|]>0$.
We deduce from the Carbery-Wright's inequality (\ref{ff1})  the existence of a universal constant $c>0$ such that, for any $n\geq 1$,
\[
P(|Q_n(I_1(e_1),\ldots,I_1(e_n)| \leq \lambda)\leq cN\lambda^{1/N}(E[Q_n(I_1(e_1),\ldots,I_1(e_n)^2])^{-1/2N}.
\]
Using the property
\[
E[Q_n(I_1(e_1),\ldots,I_1(e_n) ^2]\geq (E[|Q_n(I_1(e_1),\ldots,I_1(e_n)| ])^2
\]
we obtain
\[
P(|Q_n(I_1(e_1),\ldots,I_1(e_n)| \leq \lambda)\leq cN\lambda^{1/N}(E[| Q_n(I_1(e_1),\ldots,I_1(e_n))|])^{-1/N},
\]
and  letting $n$ tend to infinity we get
\begin{equation}\label{carbery}
P(\det \Gamma \leq \lambda)\leq cN\lambda^{1/N}(E[\det\Gamma])^{-1/N}.
\end{equation}
Letting $\lambda\to 0$, we get that $P(\det \Gamma =0)=0$. As an immediate consequence of absolute continuity criterion, (see, for instance, \cite[Theorem 2.1.1]{nualartbook}) we get the absolute continuity of the law of $F$, and assertion $(a)$ does not hold.

It is worthwhile noting that, in passing, we have proved that $P(\det\Gamma =0)$  is zero or one.\\

\noindent
{\it Proof of (b)$\Rightarrow$(a)}.
Assume the existence of $H\in\R[X_1,\cdots,X_d]\setminus\{0\}$ such that, almost surely, $H(F_1,\dots,F_d)=0$. Since $H\not\equiv 0$, the zeros of $H$ constitute a closed subset of $\R^d$ with Lebesgue measure $0$. As a result, the vector $F$ cannot have a density with respect to the Lebesgue measure.\\

\noindent
{\it Proof of (c)$\Rightarrow$(b)}.
Let $\{e_k, k\geq 1\}$ be an orthonormal basis of $\mathfrak{H}$, and set $G_k=I_1(e_k)$ for any $k\geq 1$.  In order to illustrate the method of proof, we are going  to deal first with the finite dimensional case, that is, when $F_i=P_i(G_1,\ldots,G_n)$, $i=1,\ldots,d$, and for each $i$, $P_i \in \R[X_1,\ldots,X_n]$ is a polynomial of degree at most $q$.  In that case,
\[
\langle DF_i ,DF_k \rangle_{\mathfrak{H}} = \sum_{j=1}^n
\frac{\partial P_i}{\partial x_j}(G_1,\ldots,G_n)\frac{\partial P_k}{\partial x_j}(G_1,\ldots,G_n),
\]
and the Malliavin matrix $\Gamma$ of $F$ can be written as
$\Gamma =AA^T$, where
\[
A=\left(
\frac{\partial P_i}{\partial x_j}(G_1,\ldots,G_n)
\right)_{
1\leq i\leq d,\,\,
1\leq j\leq n}.
\]
As a consequence, taking into account that the support of the law of $(G_1,\ldots,G_n)$ is $\R^n$, if $\det \Gamma =0$ almost surely, then the Jacobian
$\left(
\frac{\partial P_i}{\partial x_j}(y_1,\ldots,y_n)
\right)_{
d\times n}$
has rank strictly less than $d$ for all  $(y_1,\ldots,y_n)\in\R^n$.
Statement (b) is then a consequence of Theorem 2 and Theorem 11 in \cite{kayal}.

\medskip
Consider now the general case.  Any symmetric element $f\in \mathfrak{H}^{\otimes k}$ can be written as
\[
f=\sum_{l_1,\ldots,l_k=1}^\infty a_{l_1,\ldots,l_k} \,e_{l_1}\otimes\ldots\otimes e_{l_k}.
\]
Setting $k_l=\#\{j:\,l_j=l\}$, the multiple stochastic integral  of $e_{l_1}\otimes\ldots\otimes e_{l_k}$ can be written in terms of Hermite polynomials as
\[
I_k(e_{l_1}\otimes\ldots\otimes e_{l_k})=\prod_{l=1}^\infty H_{k_l}(G_l),
\]
where the above product is finite. Thus,
\[
I_k(f)=\sum_{l_1,\ldots,l_k=1}^\infty a_{l_1,\ldots,l_k} \prod_{l=1}^\infty H_{k_l}(G_l),
\]
where the series converges in $L^2$. This implies that we can write
\begin{equation}\label{repre-poly}
I_k(f)=P(G_1,G_2,\ldots)
\end{equation}
where $P:\R^\N\to\R$  is a function defined  $\nu^{\otimes \N}$-almost everywhere, with  $\nu$ the standard normal distribution.  In other words, we can consider $I_k(f)$ as a random variable defined in the probability space $( \mathbb{R}^{\mathbb{N}}, \nu^{\otimes \N})$.
On the other hand, for any  $n\geq 1$ and for almost all $y_{n+1},y_{n+2},\ldots$ in $\R$,
the function
$(y_1,\ldots,y_{n })\mapsto P(y_1,y_2,\ldots)$  is a polynomial of degree at most $p$.
By linearity, from the representation  (\ref{repre-poly}) we deduce the existence of mappings    $P_1,\ldots,P_d:\R^\N\to\R$, defined  $\nu^{\otimes \N}$ almost everywhere, such that
for all $i=1,\dots, d$,
\begin{equation}\label{repre-poly2}
F_i=P_i(G_1,G_2,\ldots),
\end{equation}
and such that for all  $n\geq 1$ and almost all  $y_{n+1},y_{n+2},\ldots$ in $\R$,
the mapping
$(y_1,\ldots,y_{n})\mapsto P_i(y_1,y_2,\ldots)$
is a polynomial of degree at most  $q$.
With this notation, the Malliavin matrix $\Gamma $ can be expressed as $\Gamma=AA^T$, where
\[
A=\left(
\frac{\partial P_i}{\partial x_j}(G_1,G_2,\ldots)
\right)_{
1\leq i\leq d,\,\,
j\geq 1}.
\]
Consider the truncated Malliavin matrix $\Gamma_n= A_n A_n^T$, where
\[
A_n=\left(
\frac{\partial P_i}{\partial x_j}(G_1,G_2,\ldots)
\right)_{
1\leq i\leq d,\,\,
1\leq j\leq n}.
\]
From the Cauchy-Binet formula
\begin{eqnarray*}
\det \Gamma_n = \det (A_nA_n^T)=\sum_{J=\{j_1,\ldots,j_d\}\subset \{1,\ldots,n\}} (\det A_J)^2,
\end{eqnarray*}
where for  $J=\{j_1,\ldots,j_d\}$,
\[
A_J=\left(
\frac{\partial P_i}{\partial x_j}(G_1, G_2, \ldots)
\right)_{
1\leq i\leq d,\,\,
j\in J},
\]
we deduce that $\det \Gamma_n$ is increasing and it converges to $\det \Gamma$.
Therefore, if $\det \Gamma =0$ almost surely, then for each $n\geq 1$, $\det \Gamma_n=0$ almost surely.

Suppose that   $E[\det \Gamma]=0$, which implies that $\det\Gamma =0$ almost surely. Then, for all $n\geq 1$, $\det \Gamma_n=0$ almost surely.
We can assume that  for any subset $\{F_{i_1}, \dots, F_{i_r}\}$ of the random variables $\{F_1, \dots, F_d\}$ we have
\[
E[\det \Gamma_n(F_{i_1}, \dots, F_{i_r})]\not=0,
\]
because otherwise we will work with a proper subset of this family. This implies that for $n\geq n_0$,  and for any subset $\{F_{i_1}, \dots, F_{i_r}\}$,
\[
E[\det \Gamma_n(F_{i_1}, \dots, F_{i_r})]\not=0,
\]
where $\Gamma_n$ denotes the truncated Malliavin matrix defined above.  Then,
 applying the Carbery-Wright inequality we can show that the probability
 $P(\det \Gamma_n(F_{i_1}, \dots, F_{i_r})=0)$ is zero or one, so we deduce
 $\det \Gamma_n(F_{i_1}, \dots, F_{i_r})>0$ almost surely.

Fix $n\geq n_0$. We are going to apply the results by Kayal (see \cite{kayal}) to the family of random polynomials
\[
P_i^{(n)}(y_1,\ldots,y_n)=P_i(y_1,\ldots,y_n,G_{n+1},G_{n+2},\ldots), \quad 1\leq i \leq d.
\]
We can consider these polynomials as elements of the ring of polynomials $\mathbb{K}[y_1, \dots, y_n]$, where $\mathbb{K}$ is the field generated by all multiple stochastic integrals.
This field is well defined because by a result of Shigekawa \cite{Shi} if $F$ and $G$ are finite sums of multiple stochastic integrals and $G\not \equiv 0$, then $G$ is different from zero almost surely and $\frac FG$ is well defined. The Jacobian of this set of polynomials
\[
J(y_1, \dots, y_n)=\left(
\frac{\partial P^{(n)}_i}{\partial y_j}(y_1, \dots, y_n)
\right)_{
1\leq i\leq d,\,\,
1\leq j \leq n}
\]
satisfies $J(G_1, \dots, G_n)= A_n$ almost surely, and, therefore, it has determinant zero almost surely. Furthermore, for any proper subfamily of polynomials $\{P^{(n)}_{i_1},
\dots, P^{(n)}_{i_r}\}$, the corresponding Jacobian has nonzero determinant. As a consequence of the results by Kayal, there exists a nonzero irreducible  polynomial $H_n \in \mathbb{F}[x_1, \dots, x_d]$ of degree at most $D:=dq^{d-1}$, which satisfies the following properties:
 \begin{enumerate}
 \item[$(i)$] The coefficients of $H_n$ are random variables measurable with respect to the $\sigma$-field $\sigma\{G_{n+1},G_{n+2},\ldots\}$.
 \item[$(ii)$] The coefficient of the largest monomial in antilexicographic order occurring in $H_n$ is $1$.
 \item[$(iii)$] For all $y_1, \dots, y_n \in \mathbb{R}$,
 \[
H_n(P^{(n)}_1(y_1,\ldots,y_n),\ldots,P^{(n)}_d(y_1,\ldots,y_n))=0
\]
almost surely.
 \item[$(iv)$] If $A\in \mathbb{F}[x_1, \dots, x_d]$ satisfies
\[
A(P^{(n)}_1(y_1,\ldots,y_n),\ldots,P^{(n)}_d(y_1,\ldots,y_n))=0
\]
almost surely, then $A$ is a multiple of $H_n$, almost surely.
 \end{enumerate}

 If we apply property $(iii)$ to $n+1$ and substitute $y_{n+1}$ by $G_{n+1}$ we obtain
 \[
 H_{n+1}(P^{(n+1)}_1(y_1,\ldots,y_n,G_{n+1}),\ldots,P^{(n+1)}_d(y_1,\ldots,y_n, G_{n+1}))=0.
\]
From property $(iv)$ and taking into account that for any $1\leq i \leq d$,
\[
P^{(n+1)}_i(y_1,\ldots,y_n,G_{n+1})= P_i^{(n)}(y_1, \ldots, y_n)
\] almost surely, we deduce that $H_{n+1}$ is a multiple of $H_n$ almost surely.  Using the fact that $H_{n+1}$ is irreducible and normalized we deduce that $H_n= H_{n+1}$ almost surely for any $n\geq n_0$. The coefficients of these polynomials are random variables, but, in view of condition  $(i)$, and using the $0-1$ Kolmogorov law we obtain that the coefficients are deterministic. Thus, there exists a polynomial $H\in\R[X_1,\dots,X_d]\setminus\{0\}$ of degree at most $D=dq^{d-1}$ such that  $H(F_1, \dots, F_d)=0$ almost surely.
\qed

\bigskip
The condition $E[\det \Gamma ]   >0$  can be translated into  a condition on the kernels  of the  multiple integrals appearing in the expansion of each component of the random vector $F$. Consider the following simple particular cases.

\medskip
\noindent \textbf{Example 1.}
Let $(F,G)=\big(I_1(f),I_k(g))$, with $k\geq 1$.
Let $\Gamma$ be the Malliavin matrix of $(F,G)$.
Let us compute $E[\det\Gamma]$.
Applying the duality relationship  (\ref{dual}) and the fact that $\delta(DG) =-LG=kG$, where $L$ is the Ornstein-Uhlenbeck operator, we deduce
\[
 E[\|DG\|_{\HH}^2]=E[G\delta (DG)]=kE[G^2]=kk!\|g\|_{\HH^{\otimes k}}^2,
 \]
so that
\begin{eqnarray*}
E[\det\Gamma] &=& \|f\|_{\HH}^2E[\|DG\|_{\HH}^2] - E[\langle f,DG\rangle_{\HH}^2]
=\|f\|^2_{\HH}E[\|DG\|_{\HH}^2] - k^2E[I_{k-1}(f\otimes_1 g)^2]\\
&=&kk!\big(\|f\|_{\HH}^2\|g\|_{\HH^{\otimes k}}^2 - \|f\otimes_1g \|^2_{\HH^{\otimes (k-1)}}\big).
\end{eqnarray*}
We deduce that $E[\det\Gamma]>0$ if and only if $\|f\otimes_1g \|_{\HH^{\otimes (k-1)}}<\|f\|_{\HH}\|g\|_{\HH^{\otimes k}}$.
Notice that when $k=1$ the above formula for $E[\det\Gamma]$ reduces to $E[\det\Gamma]=\det C$, where $C$ is the covariance matrix of $(F,G)$.\\

\medskip
\noindent \textbf{Example 2.}
Let $(F,G)=\big(I_2(f),I_k(g))$, with $k\geq 2$.
Let $\Gamma$ be the Malliavin matrix of $(F,G)$.
Let us compute $E[\det\Gamma]$. We
have
\begin{eqnarray*}
\|DG\|^2_{\HH} &=& k^2\sum_{r=1}^{k} (r-1)!\binom{k-1}{r-1}^2 I_{2k-2r}(g\otimes_{r}g)
= \sum_{r=1}^{k} rr!\binom{k}{r}^2 I_{2k-2r}(g\otimes_{r}g),
\end{eqnarray*}
so that
\begin{eqnarray*}
\langle DF,DG\rangle_{\HH} &=& 2k\big(
I_k(f\otimes_1 g)+(k-1)I_{k-2}(f\otimes_2 g)
\big)\\
\|DF\|^2_{\HH}&=&4\|f\|^2_{\HH^{\otimes 2}}+4I_2(f\otimes_1 f)\\
\|DG\|^2_{\HH}&=&kk!\|g\|_{\HH^{\otimes k}}^2+(k-1)kk!I_2(g\otimes_{k-1}g)+ \sum_{r=3}^{k} rr!\binom{k}{r}^2 I_{2k-2r}(g\otimes_{r}g).
\end{eqnarray*}
We deduce
\begin{eqnarray}
E[\det\Gamma] &=& E[\|DF\|^2_{\HH}\|DG\|^2_{\HH}] - E[\langle DF,DG\rangle^2_{\HH}]\notag\\
&=&4kk!\|f\|^2_{\HH^{\otimes 2}}\|g\|^2_{\HH^{\otimes k}}+8(k-1)kk!\langle f\otimes_1 f,g\otimes_{k-1}g\rangle_{\HH^{\otimes 2}}-4k^2k!\|f\,\widetilde{\otimes}_1 g\|_{\HH^{\otimes k}}^2\notag\\
&&-4k(k-1)k!\|f\otimes_2 g\|^2_{\HH^{\otimes (k-2)}}\notag\\
&=&4kk!\|f\|^2_{\HH^{\otimes 2}}\|g\|^2_{\HH^{\otimes k}}+8(k-1)kk!\| f\otimes_1 g\|_{\HH^{\otimes k}}^2-4k^2k!\|f\,\widetilde{\otimes}_1 g\|^2_{\HH^{\otimes k}}\notag\\
&&-4k(k-1)k!\|f\otimes_2 g\|^2_{\HH^{\otimes (k-2)}}. \label{35}
\end{eqnarray}
Therefore, $E[\det\Gamma]>0$ if and only if the right hand side of (\ref{35}) is strictly positive.

\medskip
Consider the particular case $k=2$, that is,  $F=(I_2(f),I_2(g))$ and let $C$ be the covariance matrix of $F$. By specializing (\ref{35})   to $k=2$, we get that
\begin{equation}\label{ineq=2}
E[\det\Gamma] = 16\big(\|f\|^2_{\HH^{\otimes 2}}\|g\|^2_{\HH^{\otimes 2}}-\langle f, g\rangle^2_{\HH^{\otimes 2}}\big)+32\big(\| f\otimes_1 g\|_{\HH^{\otimes 2}}^2-\|f\,\widetilde{\otimes}_1 g\|^2_{\HH^{\otimes 2}}\big)
\geq 4\det C .
\end{equation}

We deduce an interesting result, that generalizes a well-known criterion for Gaussian pairs.
\begin{prop}\label{pair}
Let $F=(I_2(f),I_2(g))$ and let $C$ be the covariance matrix of $F$.
Then, the law of $F$ has a density if and only if $\det C>0$.
\end{prop}
{\it Proof}. If $\det C>0$ then $E[\det\Gamma]>0$ by (\ref{ineq=2}); we deduce from Theorem
\ref{ac-vectormultiple} that the law of $F$ has a density.
Conversely, if $\det C=0$ then $I_2(f)$ and $I_2(g)$ are proportional; this prevents $F$ to have a density.
\qed

\section{Convergence in law and total variation distance}

In this section we first prove an inequality between the total variation distance and the Fortet-Mourier distance for vectors in a finite sum of Wiener chaoses.

\begin{thm}\label{tv-thm}
Fix $q,d\geq 2$, and let $F_n=(F_{1,n},\ldots,F_{d,n})$ be a sequence such that $F_{i,n}\in\bigoplus_{k=1}^q \mathcal{H}_k$ for any $i=1,\ldots,d$ and $n\geq 1$. Let $\Gamma_n:=\Gamma(F_n)$ be the Malliavin matrix of $F_n$.
Assume that $F_n\overset{\rm law}{\to} F_\infty$ as $n\to\infty$ and that there exists $\beta>0$ such that $E[\det\Gamma_n]\geq \beta$ for all $n$.
Then $F_\infty$ has a density and, for any $\gamma<\frac 1{(d+1)(4d(q-1)+3)+1}$, there exists $c>0$ such that
\begin{equation}\label{ineg-fin}
d_{TV}(F_n,F_\infty)\leq c\, d_{FM}(F_n,F_\infty)^\gamma.
\end{equation}
In particular, $F_n\to F_\infty$ in total variation as $n\to\infty$.
\end{thm}

\noindent
{\it Proof}.
The proof is divided into several steps.

\medskip

{\it Step 1}. Since
$F_{i,n}\overset{\rm law}{\to} F_{i,\infty}$ with $F_{i,n}\in\bigoplus_{k=1}^q \mathcal{H}_k$,
it follows from \cite[Lemma 2.4]{NouPol} that for any $i=1,\ldots,d$, the sequence
$(F_{i,n})$ satisfies $\sup_n E|F_{i,n}|^p<\infty$ for all $p\geq 1$. Let $\phi:\R^d\to\R\in  \mathcal{C}^\infty$ be such that $\|\phi\|_\infty\leq 1$. We can write, for any $n,m,p,M\geq 1$,
\begin{eqnarray*}
\big| E[\phi(F_n)]-E[\phi(F_m)]\big|&\leq&
\left|
E\left[
(\phi{\bf 1}_{[-M/2,M/2]^d})(F_n)
\right]
-E\left[
(\phi{\bf 1}_{[-M/2,M/2]^d})(F_m)
\right]
\right|\\
&&+2\,\sup_{n\geq 1} P(\max_{1\leq i\leq d}|F_{i,n}|\geq M/2)\\
&\leq&
\sup_{
\stackrel{\psi\in\mathcal{C}^\infty:\,\|\psi\|_\infty\leq 1}{{\rm supp}\psi\subset [-M,M]^d}
}\big| E[\psi(F_n)]-E[\psi(F_m)]\big|\\
&&+\frac{2^{1+p}}{M^p}\,\sup_{n\geq 1} E\left[\max_{1\leq i\leq d}|F_{i,n}|^p\right].
\end{eqnarray*}
Therefore, since $\sup_{n\geq 1} E\left[\max_{1\leq i\leq d}|F_{i,n}|^p\right]$ is finite,
there exists a constant $c>0$ (depending on $p$) satisfying, for all $n\geq 1$,
\begin{equation}\label{limit334}
d_{TV}(F_n,F_\infty)\leq
\sup_{
\stackrel{\phi\in\mathcal{C}^\infty:\,\|\phi\|_\infty\leq 1}{{\rm supp}\phi\subset [-M,M]^d}
}\big| E[\phi(F_n)]-E[\phi(F_\infty)]\big|
+\frac{c}{M^p}.
\end{equation}
As in \cite{NouPol}, now the idea to bound the first term in the right-hand side of (\ref{limit334}) is to regularize the function $\phi$ by means of an approximation of the identity and then to control the error term using the integration by parts of Malliavin calculus.
Let $\phi:\R^d\to\R\in\mathcal{C}^\infty$ with compact support in $[-M,M]^d$ and satisfying $\|\phi\|_\infty\leq 1$.
Let $n,m\geq 1$ be integers.
Let $0<\alpha\leq 1$ and let $\rho:\R^d\to\R_+$ be in $\mathcal{C}^\infty_c$ and satisfying $\int_{\R^d}\rho(x)dx=1$.
Set $\rho_\alpha(x)=\frac1{\alpha^d}\rho(\frac{x}{\alpha})$.
By \cite[(3.26)]{NouPol}, we have that
$\phi*\rho_\alpha$ is bounded by 1 and is Lipschitz continuous with constant $1/\alpha$.
We can thus write,
\begin{eqnarray}
&&\big|
E[\phi(F_n)]-E[\phi(F_m)]
\big|\notag\\
&\leq&\left|
E\left[\phi*\rho_\alpha(F_n)-\phi*\rho_\alpha(F_m)\right]
\right|  +2 \sup_{n\geq 1}
\left|
E\left[\left(\phi(F_n)-\phi*\rho_\alpha(F_n)\right) \right] \right| \notag\\
&\leq&\frac{1}{\alpha}d_{FM}(F_n,F_m)+
2  R_\alpha,  \label{9}
\end{eqnarray}
where $d_{FM}$  is the Forter-Mourier distance and
\[
R_\alpha= \sup_{n\geq 1}
\left|
E\left[\left(\phi(F_n)-\phi*\rho_\alpha(F_n)\right) \right] \right|.
\]
In order to estimate the term $R_\alpha$ we decompose the expectation into two parts using the identity
\[
1=\frac{\e}{\det\Gamma_n +\e}+\frac{\det\Gamma_n }{\det\Gamma_n +\e},\quad\e>0.
\]

\medskip

 {\it Step 2}. We claim that there exists $c>0$ such that, for all $\e>0$ and all $n\geq 1$,
\begin{equation}\label{cw-ineq1}
E\left[\frac{\e}{\det\Gamma_n+\e}\right]\leq c\,\e^{\frac{1}{2(q-1)d+1}}.
\end{equation}
Indeed, for any $\lambda>0$
and by using (\ref{carbery}) together with the assumption
$E[\det\Gamma_n]\geq \beta$,
\begin{eqnarray*}
E\left[\frac{\e}{\det\Gamma_n+\e}\right]\leq
E\left[\frac{\e}{\det\Gamma_n +\e}\,{\bf 1}_{\{\det\Gamma_n>\lambda\}}\right]
+c\,\lambda^{\frac{1}{2(q-1)d}}
\leq \frac{\e}{\lambda}+c\,\lambda^{\frac{1}{2(q-1)d}}.
\end{eqnarray*}
Choosing $\lambda=\e^{\frac{2(q-1)d}{2(q-1)d+1}}$
proves the claim (\ref{cw-ineq1}). As a consequence, the estimate (\ref{cw-ineq1}) implies
\begin{eqnarray}
R_\alpha&=&\sup_{n\geq 1} \left|
E\left[\left(\phi(F_n)-\phi*\rho_\alpha(F_n)\right)\left(\frac{\e}{\det\Gamma_n+\e}+\frac{\det\Gamma_n}{\det\Gamma_n +\e}\right)\right]
\right|\notag\\
&\leq& 2c\,\e^{\frac{1}{2(q-1)d+1}}+
\sup_{n\geq 1}\left|
E\left[(\phi-\phi*\rho_\alpha)(F_n)\,\frac{\det\Gamma_n}{\det\Gamma_n+\e}\right]
\right|.  \label{7}
\end{eqnarray}

\medskip

{\it Step 3}.
In this step we will derive the integration by parts formula that will be useful for our purposes.
The method is based on the representation of the density of a Wiener functional using the Poisson kernel obtained by Malliavin and Thalmaier in \cite{MaTha}, and it has been further developed by Bally and Caramellino in the works \cite{BaCa1} and \cite{BaCa2}.

Let $h:\R^d\to\R$ be a function in $\mathcal{C}^\infty$ with compact support, and consider a random variable $W\in \mathbb{D}^\infty$. Consider the Poisson kernel in $\mathbb{R}^d$ ($d\geq 2$), defined as the solution to the equation $\Delta Q_d=\delta_0$. We know that
$Q_2(x)= c_2 \log |x|$ and that  $Q_d(x) =c_d |x|^{2-d}$ for $d\geq 3$.
Then, we have the following identity
\begin{equation} \label{bb1}
h= \sum_{i=1}^d \partial_i h * \partial _i Q_d.
\end{equation}
As a consequence, we can write
\begin{eqnarray*}
E[W\det \Gamma_n\,  h(F_n)] &=&
\sum_{i=1}^d E\left[ W \det\Gamma_n\int_{\mathbb{R}^d} \partial _i Q_d(y) \partial_i h(F_n-y)dy\right]\\
&=& \sum_{i=1}^d \int_{\mathbb{R}^d} \partial _i Q_d(y)  E\left( W \det\Gamma_n  \partial_i h(F_n-y)\right)dy.
\end{eqnarray*}
We claim that
\begin{equation} \label{6}
E\left[W  \det\Gamma_n  \partial_i h(F_n-y)\right]=
\sum_{a=1}^d E \left[ h(F_n-y) \delta ( W (\mathrm{Com} \Gamma_n)_{i,a} DF_{a,n}) \right],
\end{equation}
where $\delta$ is the divergence operator, and  $\mathrm{Com}(\cdot)$ stands for the usual comatrice operator.    The equality
(\ref{6}) follows easily from the relation
\[
\partial_i h(F_n-y)= \sum_{a=1}^d(\Gamma_n^{-1})_{a,i}  \langle D(h(F_n-y)), DF_{a,n} \rangle_{\HH},
\]
multiplying by  $W \det\Gamma_n$, taking the mathematical expectation, and applying the duality relationship between the derivative and the divergence operator. The random variable
\begin{eqnarray*}
A_{i,n}(W)&=&\sum_{a=1}^d\delta ( W (\mathrm{Com} \Gamma_n)_{i,a} DF_{a,n})\\
&=&-\sum_{a=1}^{d}\big(
\langle D(W({\rm Com}\Gamma_n)_{a,i}),DF_{a,n}\rangle_\HH+({\rm Com}\Gamma_n)_{a,i}WLF_{a,n}
\big)
\end{eqnarray*}
satisfies $A_{i,n}(W)\in \mathbb{D}^\infty$, and we can write
\begin{equation}  \label{eq1}
E[W\det \Gamma_n h(F_n)]= \sum_{i=1}^d E\left[ A_{i,n}(W)\int_{\mathbb{R}^d} h(y) \partial _i Q_d(F_n-y)dy \right].
\end{equation}

\medskip

{\it Step 4}.
We are going to apply the identity (\ref{eq1}) to the function $h=\phi-\phi*\rho_\alpha$ and to the random variable $W=W_{n,\e}=\frac{ 1}{\det\Gamma_n +\e}$. In this way we obtain
\begin{eqnarray}
 &&E\left[(\phi-\phi*\rho_\alpha)(F_n)\frac{\det\Gamma_n}{\det\Gamma_n +\e}\right]\notag\\
&=&\sum_{i=1}^d E\left[ A_{i,n}(W_{n,\e})\int_{\mathbb{R}^d} (\phi-\phi*\rho_\alpha)(y) \partial _i Q_d(F_n-y)dy \right].
\label{claim2}
\end{eqnarray}
We claim that, for any $p\geq 1$, there exists a constant $c>0$ such that
\[\sup_{n} E[|A_{i,n}(W_{n,\e})|^p] \leq c \e^{-2}.
\]
Indeed, this follows immediately from the fact that the sequence $(F_{i,n})$ is uniformly bounded in $L^p$  for each $i=1,\dots,d$ and that  we can write
\begin{eqnarray*}
A_{i,n}(W_{n,\e})&=&\sum_{a=1}^{d}\Bigg\{ -\frac{ 1}{\det\Gamma_n +\e}
\left(
\langle D({\rm Com}\Gamma_n)_{a,i},DF_{a,n}\rangle_\HH-({\rm Com}\Gamma_n)_{a,i}     LF_{a,n}\right)\\
&&+\frac {1} { ( \det\Gamma_n +\e )^2}  ({\rm Com}\Gamma_n)_{a,i}
\langle D(\det\Gamma_n),DF_{a,n}\rangle_\HH  \Bigg\}.
\end{eqnarray*}
On the other hand,  we have
\begin{eqnarray*}
\int_{\mathbb{R}^d} (\phi-\phi*\rho_\alpha)(y) \partial _i Q_d(F_n-y)dy  &=&
\int_{\mathbb{R}^{2d}} (\phi(y)-\phi(y-z))  \rho_\alpha(z) \partial _i Q_d(F_n-y)dy dz\\
&=&\int_
{\mathbb{R}^{2d}} \phi(F_n-y)  \rho_\alpha(z)( \partial _i Q_d(y)-\partial_i Q_d (y-z))dy dz.
\end{eqnarray*}
Taking into account that
 \begin{equation} \label{a1}
 \partial_i Q_d(x)= k_d \frac {x_i}{ |x|^d},
 \end{equation}
 for some constant $k_d$, we can write
 \[
 \int_{\mathbb{R}^d} (\phi-\phi*\rho_\alpha)(y) \partial _i Q_d(F_n-y)dy
 = k_d
 \int_{\mathbb{R}^{2d}} \phi(F_n-y)  \rho_\alpha(z)\left( \frac{y_i}{|y|^d} -\frac{y_i-z_i}{|y-z|^d}\right)dy dz.
\]
 Fix $R>0$. Set $B_R=\{(y,z): |y|\geq R, |y-z| \geq R\}$. We can assume that the support of $\rho$ is the unit ball  $\{|z|\leq 1\}$.
 Then for any $(y,z)\in (B_R)^c$ with $|z|\leq \alpha$,
 both $|y|$ and $|y-z|$ are bounded by $R+\alpha$, and we obtain
 \begin{eqnarray*}
 &&\left| \int_
  {(B_R)^c} \phi(F_n-y)  \rho_\alpha(z)\left( \frac{y_i}{|y|^d} -\frac{y_i-z_i}{|y-z|^d}\right)dy dz\right| \\
&\leq& \int_
  {(B_R)^c}\rho_\alpha(z)\left( \frac{|y_i|}{|y|^d} +\frac{|y_i-z_i|}{|y-z|^d}\right)dy dz  \\
  &\leq& 2\int_{\{ |y| \leq R+\alpha\}}\frac{|y_i|}{|y|^d} dy
  =2\int_{\{|y|\leq 1\}}\frac{|y_i|}{|y|^d} dy\times (R+\alpha).
\end{eqnarray*}
On the other hand,
 \begin{eqnarray*}
&& \left| \int_
{B_R} \phi(F_n-y)  \rho_\alpha(z)\left( \frac{y_i}{|y|^d} -\frac{y_i-z_i}{|y-z|^d}\right)dy dz\right|\\
&& \leq  (\max_{1\leq i\leq d}|F_{i,n}|+M)^d \sup_{|y| \geq R} \int_{\{ z: |y-z| \geq R\} }
\left| \frac{y_i}{|y|^d} -\frac{y_i-z_i}{|y-z|^d}\right| \rho_\alpha(z)dz.
\end{eqnarray*}
There exists a constant $c>0$ such that,
for $|y|\geq R$, $|y-z|\geq R$ and $|z|\leq \alpha$,
\[
\left| \frac{y_i}{|y|^d} -\frac{y_i-z_i}{|y-z|^d}\right|
\leq  \frac { |y|  | |y-z|^d -|y|^d|}{ |y|^d |y-z|^d} + \frac {|z| |y|^d} { |y|^d |y-z|^d}
\leq cR^{-d} \alpha.
\]
Therefore,
\[
\left|\int_{\mathbb{R}^d} (\phi-\phi*\rho_\alpha)(y) \partial _i Q_d(F_n-y)dy\right|\leq c \left( R+ \alpha+ \alpha R^{-d} (\max_{1\leq i\leq d}|F_{i,n}|+M)^d \right),
\]
for some constant $c>0$.
Substituting this estimate into (\ref{claim2}) and assuming that $M\geq 1$,  yields
 \[
 \sup_{n} \left|E\left[(\phi-\phi*\rho_\alpha)(F_n)\frac{\det\Gamma_n}{\det\Gamma_n +\e}\right]  \right|\leq
 c\e^{-2} \left( R+ \alpha+ \alpha R^{-d}  M^d  \right),
 \]
 for some constant $c>0$.
Choosing $R=\alpha^{\frac 1{d+1}} M^{\frac d {d+1}}$ and assuming $\alpha \leq 1$, we obtain
 \begin{equation} \label{eq2}
 \sup_{n} \left|E\left[(\phi-\phi*\rho_\alpha)(F_n)\frac{\det\Gamma_n}{\det\Gamma_n+\e}\right]  \right|\leq
c\e^{-2}  \alpha^{\frac 1{d+1}} M^{\frac {d}{d+1}},
 \end{equation}
 for some constant $c>0$.

\medskip

{\it Step 5}. From (\ref{9}),  (\ref{7}) and (\ref{eq2}) we obtain
\[
|E[\phi(F_n)]-E[\phi(F_m)] |\leq \frac{1}{\alpha}d_{FM}(F_n,F_\infty)+c\e^{\frac{1}{2(q-1)d+1}}+
c\e^{-2} \alpha^{\frac 1{d+1}}  M^{\frac {d}{d+1}}.
\]
By letting $m\to\infty$, we get
\begin{equation}\label{crucial-ineq}
|E[\phi(F_n)]-E[\phi(F_\infty)]|
\leq  \frac{1}{\alpha}d_{FM}(F_n,F_\infty)+c\e^{\frac{1}{2(q-1)d+1}}+
c\e^{-2}  \alpha^{\frac 1{d+1}}  M^{\frac {d}{d+1}}.
\end{equation}
Finally, by plugging (\ref{crucial-ineq}) into (\ref{limit334}) we obtain the following inequality, valid for every $M\geq 1$, $p\geq 1$, $n\geq 1$, $\e>0$ and $0<\alpha \leq 1$:
\begin{equation}{\label{optimize1}}
d_{TV}(F_n,F_\infty)\leq c \left(\frac{1}{\alpha}d_{FM}(F_n,F_\infty)+\,\e^{\frac{1}{2(q-1)d +1}}+
\frac{\,\alpha^{\frac 1 {d+1}}\,M^{\frac d{d+1}}}{\e^{2}}+\frac{1}{M^p}\right),
\end{equation}
where the constant $c$ depends on $p$.
Choosing  $\e= \left( \alpha^{\frac 1 {d+1}}\,M^{\frac d{d+1}} \right)^{\frac { 2(q-1)d+1}{ 2(2(q-1)d+1)+1}}$ we get
\begin{equation}{\label{optimize2}}
d_{TV}(F_n,F_\infty)\leq c \left(\frac{1}{\alpha}d_{FM}(F_n,F_\infty)+ \left( \alpha^{\frac 1 {d+1}}\,M^{\frac d{d+1}} \right)^
{\frac 1{ 4(q-1)d+3}} + M^{-p} \right).
\end{equation}
Notice that  $d_{FM}(F_n,F_\infty) \leq 1$ for $n$ large enough ($n\geq n_0$ say). So, assuming that $n\geq n_0$ and choosing
\[
\alpha=   d_{FM}(F_n,F_\infty) ^{\frac {(d+1)(4(q-1)d+3)}{(d+1)(4(q-1)d+3)+1}} M^{-\frac d  {(d+1)(4(q-1)d+3)+1}},
\]
we obtain
\begin{equation}{\label{optimize3}}
d_{TV}(F_n,F_\infty)\leq c \left(   d_{FM}(F_n,F_\infty) ^{\frac 1 D} M^{ \frac {d}D} +   M^{-p} \right),
\end{equation}
where $D= (d+1)(4(q-1)d+3)+1$. Notice that $\alpha \leq 1$ provided $M\geq 1$ and  $n\geq n_0$.
Optimizing with respect to $M$ yields
\[
d_{TV}(F_n,F_\infty)\leq c   d_{FM}(F_n,F_\infty) ^{\frac p {pD+ d}},
\]
and taking into account that $p$ can be chosen arbitrarily large, we have proved that for any $\gamma <\frac 1D$ there exists $c>0$ such that
 (\ref{ineg-fin}) holds true.

\medskip

{\it Step 6}.
Finally, let us prove that the law of $F_\infty$ is absolutely continuous with respect to the Lebesgue measure. Let $A\subset\R^d$ be a Borel set of Lebesgue measure zero.
By Lemma \ref{ac-vectormultiple} and because $E[\det\Gamma_n]\geq \beta>0$, we have $P(F_n\in A)=0$.
Since $d_{TV}(F_n,F_\infty)\to 0$ as $n\to\infty$, we deduce that $P(F_\infty\in A)=0$,
proving that $F_\infty$ has a density by
the Radon-Nikodym theorem. The proof of the theorem is now complete.

\qed

Under the assumptions of Theorem  \ref{tv-thm}, if we denote by  $\rho_n$ (resp. $\rho_\infty$) the density of $F_n$ (resp. $F_\infty$), then the convergence in total variation is equivalent to
\[
\int_{\R^d} |\rho_n(x)-\rho_\infty(x)|dx\to 0,
\]
as $n$ tends to infinity.  We are going to show that this convergence actually holds in $L^p(\mathbb{R}^d)$ for any $1\leq p< 1+\frac{1}{2d^2(q-1)+d-1}$.

\begin{prop}
Suppose that $F_n$ is a sequence of $d$-dimensional random vectors satisfying the conditions of Theorem \ref{tv-thm}. Denote by  $\rho_n$ (resp. $\rho_\infty$) the density of $F_n$ (resp. $F_\infty$). Then, for any $1\leq p<1+\frac{1}{2d^2(q-1)+d-1}$, we have
\[
\int_{\R^d} |\rho_n(x)-\rho_\infty(x)|^pdx\to 0.
\]
\end{prop}
{\it Proof}. The proof will be done in several steps. We set $N=2d(q-1)$
and we fix $p$ such that $1< p<1+\frac{1}{2d^2(q-1)+d-1}$.

1) Denote by $\Gamma_n$ the Malliavin matrix of $F_n$. Using Carbery-Wright's inequality (\ref{ff1}), we have, for any $\gamma<\frac1 N$,
\begin{eqnarray*}
\sup_n E\left[ (\det \Gamma_n )^{-\gamma}\right]=\sup_n\int_0^\infty \gamma t^{\gamma-1}P(\det \Gamma_n <t^{-1})dt
\leq C\left(1+\int_1^\infty t^{\gamma-1-\frac1{N}}dt\right)<\infty.
\end{eqnarray*}

2) Fix a real number $M>0$. For any  $\alpha<\frac 1{N+1}$ and $1+\frac \alpha d<p< 1+ \frac \alpha{d-\alpha}$, we have
\begin{eqnarray}  \notag
 &&\int_{\R^d} \rho_{n}^p(x){\bf 1}_{\{|\rho_n(x)|\leq M\}}dx\\
 \notag&=&E\left[\rho_n^{p-1}(F_n)
 {\bf 1}_{\{|\rho_n(F_n)|\leq M\}}
 \frac{(\det\Gamma_n)^{\alpha}}{(\det\Gamma_n)^{\alpha}}\right]\\ \notag
&\leq& E\left[\rho_n^{(p-1)(N+1)}(F_n){\bf 1}_{\{|\rho_n(F_n)|\leq M\}}(\det\Gamma_n)^{(N+1)\alpha}\right]^{\frac{1}{N+1}} E\left[(\det \Gamma_n)^{-\frac{N+1}{N}\alpha}\right]^{\frac{N}{N+1}}\\  \label{bb2}
&\leq&C E\left[\rho_n^{\frac{p-1}\alpha}(F_n) {\bf 1}_{\{|\rho_n(F_n)|\leq M\}}\det \Gamma_n\right]^{\alpha},
\end{eqnarray}
where
\[
C:=\sup_nE\left[(\det \Gamma_n)^{-\frac{N+1}{N}\alpha}\right]^{\frac{N}{N+1}} <\infty.
\]
Applying the identity (\ref{eq1}) to $h=\rho_n^{\frac{p-1}\alpha}{\bf 1}_{\{|\rho_n(\cdot)|\leq M\}}$ and $W=1$ and taking into account that (\ref{a1}) holds, yields
\begin{equation}
E\left[\rho_n^{\frac{p-1}\alpha}(F_n) {\bf 1}_{\{|\rho_n(F_n)|\leq M\}}\det \Gamma_n\right] \label{bb3}
= k_d \sum_{i=1}^d \int_{\mathbb{R}^d}  \frac {y_i} {|y|^d}
E\left[ \rho_n^{\frac{p-1}\alpha}(F_n-y)
{\bf 1}_{\{|\rho_n(F_n-y)|\leq M\}}A_{i,n}  \right]dy,
\end{equation}
where  $A_{i,n}=A_{i,n}(1)= \sum_{a=1}^d \delta\left( (\mathrm{Com} \Gamma_n)_{i,a} DF_{a,n} \right) $.

For any $x\in \mathbb{R}^d$ and any function $f:\mathbb{R}^d\rightarrow \mathbb{R}_+$, the integral
\[
 \int_{\R^d}  \frac {y_i} {|y|^d}  f(x-y) dy
 \]
 can be decomposed into the regions $\{y: |y| \leq 1\}$ and $\{ y: |y| >1\}$. Then, using H\"older's inequality, for any exponents $\beta>d$ and $\gamma <d$, there exist a constant $C_{\beta, \gamma}$ such that
 \[
\sup_{x\in \mathbb{R}^d} \left| \int_{y\in\R^d}  \frac {y_i} {|y|^d}  f(x-y) dy \right|\le
C_{\beta, \gamma} \left(  \|f\|_\beta+  \|f\|_\gamma \right).
 \]
 We are going to apply this estimate to the function  $f=\rho_n^{\frac{p-1}\alpha}{\bf 1}_{\{|\rho_n(\cdot)|\leq M\}}$ and to the exponents $\beta=\frac{p\alpha}{p-1}>d$ and $\gamma= \frac \alpha {p-1} <d$. In this way we obtain from (\ref{bb3})
 \begin{eqnarray}
  &&E\left[\rho_n^{\frac{p-1}\alpha}(F_n) {\bf 1}_{\{|\rho_n(F_n)|\leq M\}}\det \Gamma_n\right]\notag\\
  &\leq& c_d \sum_{i=1}^n E[|A_{i,n}|] C_{\alpha,\beta}\left[ \left(\int_{\mathbb{R}^d} \rho_n ^{p}(x) {\bf 1}_{\{|\rho_n(x)|\leq M\}}dx\right) ^{\frac{p-1}{p\alpha}} +\left(\int_{\mathbb{R}^d} \rho_n(x) {\bf 1}_{\{|\rho_n(x)|\leq M\}} dx\right) ^{\frac {p-1}\alpha}\right]\notag\\
  &\leq& c_d \sum_{i=1}^n E[|A_{i,n}|] C_{\alpha,\beta}\left[ \left(\int_{\mathbb{R}^d} \rho_n^{p}(x) {\bf 1}_{\{|\rho_n(x)|\leq M\}} dx\right) ^{\frac{p-1}{p\alpha}} +1\right]. \label{bb4}
\end{eqnarray}
 From (\ref{bb2})  and (\ref{bb4}) we deduce the existence of a constant $K$, independent of $M$ and $n$, such that
 \begin{equation}\label{123}
   \int_{\R^d} \rho_{n}^p(x){\bf 1}_{\{|\rho_n(x)|\leq M\}}dx \leq K\left[ \left(\int_{\mathbb{R}^d} \rho_n ^{p}(x) {\bf 1}_{\{|\rho_n(x)|\leq M\}}dx\right) ^{\frac{p-1}p}  +1\right].
   \end{equation}
Since $  \int_{\R^d} \rho_{n}^p(x){\bf 1}_{\{|\rho_n(x)|\leq M\}}dx$ is finite (it is less than $M^{p-1}$), we deduce from (\ref{123}) that
\[
\sup_n\sup_{M>0}  \,\,  \int_{\R^d} \rho_{n}^p(x){\bf 1}_{\{|\rho_n(x)|\leq M\}}dx<\infty,
\]
implying in turn that
\[
\sup_n  \int_{\R^d} \rho_{n}^p(x)dx<\infty,
\]

3) Let $n,m\geq 1$. By applying H\"older to
\[
\int_{\R^d} |\rho_n(x)-\rho_m(x)|^pdx
=
\int_{\R^d} |\rho_n(x)-\rho_m(x)|^{\epsilon}|\rho_n(x)-\rho_m(x)|^{p-\epsilon}dx,
\]
we obtain, for any $0<\epsilon<1$,
\begin{equation}\label{interpol}
\int_{\R^d} |\rho_n(x)-\rho_m(x)|^pdx
\leq \left(\int_{\R^d} |\rho_n(x)-\rho_m(x)|dx\right)^{\epsilon}
\left(\int_{\R^d} |\rho_n(x)-\rho_m(x)|^\frac{p-\epsilon}{1-\epsilon}dx\right)^{1-\epsilon}.
\end{equation}
We can choose $\epsilon>0$ small enough such that
\[
p< \frac{p-\epsilon}{1-\epsilon} < 1+\frac 1{2d^2(q-1)+d+1}.
\]
Then, from Part 2) we deduce
\[
\int_{\R^d} |\rho_n(x)-\rho_m(x)|^pdx \to 0\quad \mbox{as $n,m\to\infty$}.
\]
As a result, $\{\rho_n\}$ converges in $L^p(\R^d)$,
 which is the desired conclusion.
\qed

\section{Some applications}\label{app}

In this section we present some consequences of Theorems \ref{ac-vectormultiple} and \ref{tv-thm}.
We start with a straightforward consequence of Theorem \ref{tv-thm}.

\begin{prop}\label{azerty}
Fix $q,d\geq 2$, and let $F_n=(F_{1,n},\ldots,F_{d,n})$ be a sequence such that $F_{i,n}\in\bigoplus_{k=1}^q \mathcal{H}_k$ for any $i=1,\ldots,d$ and $n\geq 1$. Let $\Gamma_n:=\Gamma(F_n)$ be the Malliavin matrix of $F_n$. As $n\to\infty$, assume that $F_n\to F_\infty$ in law and that
$\Gamma_n\to M_\infty$ in law, with $E[\det M_\infty]>0$.
Then $F_n$ converges to $F_\infty$ in total variation.
\end{prop}
{\it Proof}.
We set $N=2d(q-1)$.
Since
$\Gamma_{n}(i,j)\overset{\rm law}{\to} M_{\infty}(i,j)$ with $\Gamma_{n}(i,j)\in\bigoplus_{k=0}^N \mathcal{H}_k$,
it follows from \cite[Lemma 2.4]{NouPol} that for any $i,j=1,\ldots,d$, the sequence
$\Gamma_{n}(i,j)$, $n\geq 1$, satisfies $\sup_n E|\Gamma_{n}(i,j)|^p<\infty$ for all $p\geq 1$.
As a result, $E[\det \Gamma_n]\to E[\det M_\infty]>0$ and the desired conclusion follows from Theorem \ref{tv-thm}.
\qed

\bigskip

Our first application of Proposition \ref{azerty} consists in strengthening the celebrated Peccati-Tudor \cite{PTu04} criterion of asymptotic normality.
\begin{thm}\label{PT}
Let $d\geq 2$ and $k_1, \ldots, k_d\geq 1$ be some fixed
integers. Consider vectors
\[
F_n=(F_{1,n},\ldots,F_{d,n})=
(I_{k_1}(f_{1,n}),\ldots,I_{k_d}(f_{d,n})), \quad n\geq 1,
\]
with
$f_{i,n}\in \EuFrak{H}^{\odot k_i}$.
As $n\to\infty$, assume that $F_n \overset{\rm law}{\to} N\sim\mathcal{N}_d(0,C)$ with $\det(C)>0$.
Then, $d_{TV}(F_n,N)\to 0$ as $n\to\infty$.
\end{thm}
{\it Proof}.
Since
$F_{i,n}\overset{\rm law}{\to} F_{i,\infty}$ with $F_{i,n}\in\mathcal{H}_{k_i}$,
it follows from \cite[Lemma 2.4]{NouPol} that for any $i=1,\ldots,d$, the sequence
$(F_{i,n})$ satisfies $\sup_n E|F_{i,n}|^p<\infty$ for all $p\geq 1$.
In particular, one has that $E[F_{i,n}F_{j,n}]\to C(i,j)$ as $n\to\infty$ for any
$i,j=1,\ldots,d$. Denote
by $\Gamma_n$ the Malliavin matrix of $F_n$.
As a consequence of the main result in Nualart and Ortiz-Latorre \cite{NOL}, we deduce that $\Gamma_n\to C$ in $L^2(\Omega)$ as $n\to\infty$.
Finally, the desired conclusion follows from Proposition \ref{azerty}.
\qed

\bigskip

The next result is a corollary of Proposition \ref{azerty} and Theorem \ref{ac-vectormultiple},  and improves substantially  Theorem 4 of Breton  \cite{Bre}. It represents a multidimensional version of a result by Davydov and Martynova \cite{DavMar}.

\begin{cor}
Fix $d\geq 2$ and  $k_1,\ldots,k_d\geq 1$. Consider a sequence of $d$-dimensional random vectors $\{F_n, n\geq 1\}$ of the form
 $F_n=(F_{1,n},\ldots,F_{d,n})$, with $F_{i,n}=I_{k_i}(f_{i,n})$, $i=1,\ldots,d$, $n\geq 1$.
Suppose that $F_n$ converges in $L^2(\Omega)$ to $F_\infty$ and the law of  $F_\infty$ is absolutely continuous with respect to the Lebesgue measure.  Then,  $F_n$ converges to $F_\infty$ in total variation.
\end{cor}
{\it Proof}.
By the isometry of multiple stochastic integrals, for any $i=1,\dots, d$, the sequence $f_{i,n}\in \mathfrak{H}^{\odot k_i}$ converges as $n$ tends to infinity  to  an element $f_{i,\infty} \in\mathfrak{H}^{\odot  k_i}$, and we can write
\[
F_\infty=(I_{k_1}(f_{1,\infty}),\ldots,I_{k_d}(f_{d,\infty})).
\]
Since
the law of  $F_\infty$ is absolutely continuous with respect to the Lebesgue measure,  we deduce from Theorem  \ref{ac-vectormultiple} that $E[\det\Gamma(F_\infty)]>0$, where $\Gamma(F_\infty)$ is the Malliavin matrix of $F_\infty$.  On the other hand, taking into account that all the norms
$\|\cdot \|_{m,p}$ are equivalent in a fixed Wiener chaos, we deduce that for all $1\leq i,j \leq d$
\[
\Gamma_{i,j}(F_n) \rightarrow \Gamma_{i,j}(F_\infty)
\]
 in $L^p(\Omega)$ as $n$ tends to infinity, for all $p\geq 2$. Therefore, we can conclude the proof using Proposition \ref{azerty}.\qed

\bigskip

In the case of a sequence of  $2$-dimensional vectors in the second chaos, it suffices to assume that the covariance of the limit is non singular. In fact, we have the following result.
\begin{cor}
Let $(F_n,G_n)=(I_2(f_n),I_2(g_n))$ be a pair converging in law to $(F_\infty,G_\infty)$ as $n$
tends to $\infty$.
Let $C_\infty$ be the covariance matrix of $(F_\infty,G_\infty)$ and assume that $\det C_\infty>0$.
Then $(F_n,G_n)$ converges to $(F_\infty,G_\infty)$ in total variation.
\end{cor}
{\it Proof}.
Let $\Gamma_n$ (resp. $C_n$) be the Malliavin (resp. covariance) matrix of $(F_n,G_n)$.
Taking into account that all $p$-norms are equivalent in a fixed Wiener chaos, we deduce that that both $\{F_n, n\geq 1\}$ and $\{G_n, n\geq 1\}$ are
uniformly bounded with respect to $n$ in all the $L^p(\Omega)$. Thus, one has $\det C_n\to\det C_\infty$
as $n$ tends to $\infty$.
On the other hand, we have by (\ref{ineq=2})  that $E[\det\Gamma_n]\geq 4\det C_n$.
By letting $n$ tend to $\infty$, we deduce that $E[\det\Gamma_n]\geq \frac12\det C_\infty>0$
for $n$ large enough. Theorem \ref{tv-thm} allows us to conclude.\qed

\bigskip

Another situation where we only need the limit to be non degenerate in order to obtain the convergence in total variation, is the case where the limit has pairwise independent components.

\begin{cor}\label{cor2}
Fix $d\geq 2$ and $k_1, \dots, k_d \geq 1$. Consider a sequence of $d$-dimensional random vectors of multiple stochastic integrals $\{F_n, f\geq 1\}$ of the form  $F_n=(F_{1,n},\ldots,F_{d,n})=(I_{k_1}(f_{1,n}),\ldots,I_{k_d}(f_{d,n}))$. Suppose that $F_n$
converges in law to $F_\infty=(F_{1,\infty},\ldots,F_{d,\infty})$. Assume moreover that ${\rm Var}(F_{j,\infty})>0$ for any $j=1,\dots, d$  and that $F_{1,\infty},\ldots,F_{j,\infty}$ are pairwise independent. Then
$F_n$ converges to $F_\infty$ in total variation.
\end{cor}
{\it Proof}.
The proof is divided into several steps.

\medskip

{\it Step 1}. We claim that there exists $\gamma>0$ such that $E[\|DF_{1,n}\|^2\ldots \|DF_{d,n}\|^2]\geq \gamma$ for all $n$ large enough. Indeed, let $j=1,\ldots,d$. Since $E[\|DF_{j,n}\|^2]\geq {\rm Var}(F_{j,n})$
and ${\rm Var}(F_{j,n})\to {\rm Var}(F_{j,\infty})$ as $n\to\infty$, we have that $E[\|DF_{j,n}\|^2]\geq \frac12{\rm Var}(F_{j,\infty})>0$
for all $n$ large enough. Using Carbery-Wright's inequality, we deduce that there exists $c>0$ such that, for all $n$ large enough and all $\lambda>0$,
\[
P(\|DF_{j,n}\|^2\leq \lambda)\leq c\lambda^{\frac{1}{2k_j-2}}.
\]
As a consequence, for $0<\alpha<\frac1{d(-1+\max_{1\leq j\leq d} k_j)}$, we can write
\begin{eqnarray*}
&&E[\|DF_{1,n}\|^{-\alpha}\ldots \|DF_{d,n}\|^{-\alpha}]\\
&=&\int_0^\infty P\left(\|DF_{1,n}\|^{\alpha}\ldots \|DF_{d,n}\|^{\alpha}\leq \frac1x\right)dx\\
&\leq& 1+\int_1^\infty \left\{ P\left(\|DF_{1,n}\|^{\alpha}\leq \frac{1}{x^{\frac1d}}\right)+\cdots+
P\left(\|DF_{d,n}\|^{\alpha}\leq \frac{1}{x^{\frac1d}}\right)\right\}dx\\
&\leq&c\left(1+
\int_1^\infty \left[
x^{-\frac{1}{\alpha d(k_1-1)}}+\cdots+x^{-\frac{1}{\alpha d(k_d-1)}}
\right]dx\right),
\end{eqnarray*}
so that $\sup_{n\geq 1} E[\|DF_{1,n}\|^{-\alpha}\ldots \|DF_{d,n}\|^{-\alpha}]<\infty$.
Combined with
\begin{eqnarray*}
E[\|DF_{1,n}\|^{2}\ldots \|DF_{d,n}\|^{2}]&\geq& E[\|DF_{1,n}\|^{\alpha}\ldots \|DF_{d,n}\|^{\alpha}]^{\frac{2}{\alpha}}\\
 &\geq& \big(E[\|DF_{1,n}\|^{-\alpha}\ldots \|DF_{d,n}\|^{-\alpha}]\big)^{-\frac{2}\alpha},
\end{eqnarray*}
this proves the claim.

\medskip

{\it Step 2}. We claim that $E[\det\Gamma(F_n)]-E[\|DF_{1,n}\|^2\ldots \|DF_{d,n}\|^2]\to 0$ as $n\to\infty$.
To prove the claim it suffices to show that, for any $1\leq i\neq j\leq d$, one has $\langle DF_{i,n},DF_{j,n}\rangle_\HH\to 0$
in all the $L^p(\Omega)$. By hypercontractivity, to prove this latter property it is enough to check that
$\langle DF_{i,n},DF_{j,n}\rangle_\HH\to 0$
in $L^2(\Omega)$.
Recall from \cite[(3.20)]{NouRos} that
\begin{eqnarray*}
{\rm Cov}(F_{i,n}^2,F_{j,n}^2)&=&k_i!k_j!\sum_{r=1}^{k_i\wedge k_j}\binom{k_i}{r}\binom{k_j}{r}\|f_{i,n}\otimes_r
f_{j,n}\|^2\\
&&+\sum_{r=1}^{k_i\wedge k_j}r!^2\binom{k_i}{r}^2\binom{k_j}{r}^2(k_i+k_j-2r)!\|f_{i,n}\widetilde{\otimes}_r f_{j,n}\|^2.
\end{eqnarray*}
Since ${\rm Cov}(F_{i,n}^2,F_{j,n}^2)\to {\rm Cov}(F_{i,\infty}^2,F_{j,\infty}^2)=0$ (recall that
$F_{i,\infty}$ and $F_{j,\infty}$ are assumed to be independent),
we deduce that $\|f_{i,n}\widetilde{\otimes}_r f_{j,n}\|^2\to 0$ for all $r=1,\ldots,k_i\wedge k_j$.
But
\[
E[\langle DF_{i,n},DF_{j,n}\rangle^2_\HH]
=k_i^2k_j^2\sum_{r=1}^{k_i\wedge k_j} (r-1)!^2\binom{k_i-1}{r-1}^2\binom{k_j-1}{r-1}^2(k_i+k_j-2r)!
\|f_{i,n}\widetilde{\otimes}_r f_{j,n}\|^2,
\]
and the claim is shown.

\medskip

{\it Step 3}. By combining Steps 1 and 2, we obtain the existence of $\gamma>0$ such that $E[\det\Gamma(F_n)]\geq\gamma$
for all $n$ large enough.  Theorem \ref{tv-thm} then gives the desired conclusion of Corollary \ref{cor2}.
\qed

\bigskip

\noindent
{\bf Acknowledgments}. We  thank two anonymous referees for their very careful reading of the manuscript. The first author thanks Jean-Christophe Breton for useful conversations.

\end{document}